\documentclass[11pt,leqno]{article}

\usepackage{amssymb,amsmath,amsthm,rotating}

\numberwithin{equation}{section}

\newcommand{\n}{\noindent}
\newcommand{\bb}[1]{\mathbb{#1}}
\newcommand{\cl}[1]{\mathcal{#1}}

\newcommand{\ds}{\displaystyle}

\newcommand{\ovl}{\overline}

\theoremstyle{plain}
\newtheorem{pro}{Proposition}
\newtheorem{thm}{Theorem}
\newtheorem{lem}{Lemma}

\theoremstyle{definition}

\newtheorem*{rem}{Remark}

\begin{document}

\title{On the Structure of the Square of a $C_0(1)$ 
Operator\thanks{\noindent 2000 AMS Classification:\ 47A15, 47A45.\newline
\indent ~~Keywords:\ $C_0$ operators, invariant subspace lattice.}}

\author{Ronald G.\ Douglas and Ciprian Foias\\
\em Dedicated to I.B.~Simonenko on his seventieth birthday}

\date{}
\maketitle

\setcounter{section}{-1}

\section{}

\indent

While the model theory for contraction operators (cf.\ \cite{S-F}) is always a 
useful tool, it is 
particularly powerful when dealing with $C_0(1)$ operators. Recall that 
an operator $T$ on a Hilbert space $H$ is a $C_0(N)$-operator $(N=1,2\ldots)$ if 
$\|T\|\le 1$, $T^n\to 0$ and, $T^n\to 0$ (strongly) when $n\to\infty$ and 
rank$(1-T^*T)=N$. In particular, a $C_0(1)$ operator is unitarily equivalent to 
the compression of the unilateral shift operator $S$ 
on the Hardy 
space $H^2$ to a subspace $H^2\ominus mH^2$ for some inner function $m$ in 
$H^\infty$.

In this note we use the structure theory to determine when the lattices of 
invariant and hyperinvariant subspaces differ for the square $T^2$ of a 
$C_0(1)$ 
operator and the relationship of that to the reducibility of $T^2$. To 
accomplish 
this task 
we first determine very explicitly the 
characteristic operator function for $T^2$ and use the representation obtained 
to 
determine when the operator is  irreducible. While every operator $T$ in 
$C_0(1)$ is 
irreducible, it does not follow that 
$T^2$ 
is necessarily 
irreducible, that is, has no reducing subspaces. In particular, we characterize 
those $T$ in $C_0(1)$ for which 
$T^2$ is irreducible but  for which the lattices of invariant and hyperinvariant 
subspaces 
for $T^2$ are distinct.

Finally, we provide an example of an operator $X$ on a four dimensional 
Hilbert 
space for which the two lattices are distinct but $X$ is irreducible, and show 
that such an example is 
not possible on a three dimensional space. 

This work was prompted by a 
question to the first author
from Ken Dykema (Sect.~2, \cite{D}) concerning hyperinvariant subspaces in 
von~Neumann algebras. He asked whether the lattices of invariant and 
hyperinvariant subspaces for an irreducible matrix must coincide. He provides an 
example in \cite{D} on a six-dimensional Hilbert space 
showing that this
is not the case.

We assume that the reader is familiar with the concepts and notation in \cite{B} 
and \cite{S-F}.

\section{}

\indent

Let  $T\in C_0(1)$ on $H$, $\dim H\ge 2$. WLOG we can assume
\begin{equation}\label{eq1.1}
T = P_HS|_{\ds H}, \text{ where } H = H^2 \ominus mH^2, (Sh)z= zh(z) (z\in 
D, h\in H^2), m\in H^\infty, m \text{ inner.}
\end{equation}

Define
\begin{equation}\label{eq1.2}
\Theta(\lambda) = \frac12\left[\begin{matrix}
b(\lambda)&\lambda d(\lambda)\\ d(\lambda)&b(\lambda)\end{matrix}\right]\quad 
(\lambda\in D),
\end{equation}
where
\begin{subequations}
\begin{align}
\label{eq1.3a}
&~~~~~b(\lambda) = m(\sqrt\lambda) + m(-\sqrt \lambda)\qquad (\lambda\in {\bb 
D})\quad 
\text{and}\\
\label{eq1.3b}
&\left\{\begin{array}{l}
d(\lambda) = \dfrac{m(\sqrt\lambda) - m(-\sqrt\lambda)}{\sqrt\lambda}\qquad 
(0\ne \lambda\in {\bb D})\\
d(0) = 2m'(0).\end{array}\right.
\end{align}
\end{subequations}

\begin{lem}\label{lem1}
The matrix function $\Theta(\cdot)$ is inner, pure and (up to a coincidence) 
the 
characteristic operator function of $T^2$.
\end{lem}

\begin{proof}
For $h\in H^2$ write
\begin{equation}
\label{eq1.4a}
h(\lambda) = h_0(\lambda^2) + \lambda h_1(\lambda^2)\qquad (\lambda\in {\bb 
D}).
\tag*{(1.4{\rm a})}
\end{equation}
Clearly $h_0(\cdot), h_1(\cdot)$ $(= h_0(\lambda), h_1(\lambda), \lambda\in 
{\bb 
D})$ belong to $H^2$. Define $W\colon \ H^2\mapsto H^2\oplus H^2$ $(= H^2({\bb 
C}^2))$ by
\begin{equation}\label{eq1.4b}
Wh = h_0\oplus h_1, \text{ where $h$ is given by \ref{eq1.4a}.}\tag*{(1.4{\rm 
b})}
\end{equation}
\addtocounter{equation}{1}
Then $W$ is unitary and
\begin{equation}\label{eq1.5}
WS^2 = (S\oplus S) W.
\end{equation}
Consequently,
\begin{equation}\label{eq1.6}
WT^2 = WP_HS^2 = P_{WH}WS^2 = P_{WH}(S\oplus S)W;
\end{equation}
moreover, since $S^2mH^2 \subset mH^2$ we also have
\[
(S\oplus S) WmH^2 = WS^2mH^2 \subset WmH^2
\]
and therefore
\begin{equation}\label{eq1.7}
\left\{\begin{array}{l}
P_{WH}(S\oplus S) = P_{WH}(S\oplus S) P_{WH} = WP_H W^*(S\oplus S)P_{WH}=\\
\phantom{P_{WH}(S\oplus S)} = WP_H S^2W^*P_{WH} = WT^2 P_HW^*=\\
\phantom{P_{WH}(S\oplus S)} = W|_{\ds H}T^2(W|_{\ds H})^*.\end{array}\right.
\end{equation}
These relationships show that $S\oplus S$ is an isometric lifting of 
$T_0=P_{WH}(S\oplus S)|_{\textstyle WH}$ and that this operator is unitarily 
equivalent to 
$T^2$. Moreover, since
\[
\bigvee^\infty_{n=0} (S\oplus S)^n WH = H^2\oplus H^2
\]
is obvious, $S\oplus S$ is {\em the\/} minimal isometric lifting of 
$T=W|_{\textstyle H}T^2(W|_{\textstyle H})^*$.

Further,
\begin{align*}
WmH^2 &= \{W(m_0(\lambda^2) + \lambda m_1(\lambda^2))(h_0(\lambda^2) + \lambda 
h_1(\lambda^2))\colon \ h\in H^2\}\\
&= \{W[(m_0h_0)(\lambda^2) + \lambda^2(m_1h_1)(\lambda^2) +\\
&\qquad + \lambda(m_0h_1 + m_1h_0)(\lambda^2)\colon \ h\in H^2\}=\\
&= \{((m_0h_0)(\lambda) + \lambda(m_1h_1)(\lambda)) \oplus (m_0h_1 + 
m_1h_0)(\lambda)\colon \ h\in H^2\}=\\
&= \left\{\begin{bmatrix} m_0&\lambda m_1\\ m_1&m_0\end{bmatrix} (h_0\oplus 
h_1)\colon \ h\in H^2\right\} = \begin{bmatrix} m_0&\lambda m_1\\ 
m_1&m_0\end{bmatrix} H^2\oplus H^2.
\end{align*}
Note that the above computations also prove that
\begin{equation}\label{eq1.8}
(Wm(S)W^*)(h_0\oplus h_1) = \begin{bmatrix} m_0&\lambda m_1\\ 
m_1&m_0\end{bmatrix} h_0\oplus h_1\qquad (h_0\oplus h_1\in H^2\oplus H^2).
\end{equation}
Since $m(S)$ is isometric, so is $Wm(S)W^*$, that is,
\begin{equation}\label{eq1.9}
M(\lambda) \equiv \begin{bmatrix} m_0(\lambda)&\lambda m_1(\lambda)\\ 
m_1(\lambda)&m_0(\lambda)\end{bmatrix} \text{ is inner.}
\end{equation}
Consequently, $T_0$ is the compression of $S\oplus S$ to
\begin{equation}\label{eq1.10}
WH = (H^2\oplus H^2) \ominus M(H^2\oplus H^2).
\end{equation}
Moreover, it is clear that 
\[
m_0(\lambda) = \frac12 b(\lambda),\quad m_1(\lambda) = \frac12 d(\lambda) 
\qquad 
(\lambda\in {\bb D})
\]
so that the matrix $M(\cdot)$ defined by \eqref{eq1.9} is identical to the 
matrix $\Theta(\cdot)$ defined by \eqref{eq1.2}.

Note that
\[
\Theta(0) = \begin{bmatrix} m(0)&0\\ m'(0)&m(0)\end{bmatrix}
\]
and
\[
\Theta(0)^* \Theta(0) = \begin{bmatrix} |m(0)|^2 + |m'(0)|^2&\ovl{m'(0)}m(0)\\ 
\ovl{m(0)} m'(0)&|m(0)|^2\end{bmatrix}.
\]
If $\Theta(0)$ were not pure, then $\Theta(0)^* \Theta(0)$ would have the 
eigenvalue 1 and therefore the other eigenvalue must be $|m(0)|^4$. Taking 
traces 
we have
\[
2|m(0)|^2 + |m'(0)|^2 = 1+|m(0)|^4.
\]
This implies that the modulus of the analytic function $\widetilde m(\lambda)$ 
defined by
\[
\lambda\widetilde m(\lambda) = \frac{m(\lambda)-m(0)}{1-\ovl{m(0)} m(\lambda)} 
\qquad 
(\lambda\in\ovl{\bb D}, \lambda\ne 0)
\]
and 
\[
\widetilde m(0) = \frac{m'(0)}{1-|m(0)|^2}
\]
attains its maximum $(=1)$ at $\lambda=0$. By virtue of the maximum principle, 
$\widetilde m(\lambda) = c=$ constant, $|c|=1$. Thus
\[
m(\lambda) \equiv c\left(\frac{\lambda+\bar cm(0)}{1+\lambda 
c\ovl{m(0)}}\right)\qquad (\lambda\in\ovl{\bb D})
\]
and 
\[
2\le \dim H = \dim(H^2\ominus mH^2)=1, 
\]
which is a contradiction.

We conclude that $\Theta(\cdot)$ is pure and, by virtue of \eqref{eq1.10} 
(recall $\Theta(\lambda) \equiv M(\lambda)$), that $\Theta(\cdot)$ is the 
characteristic operator function of $T_0$ and hence (up to a coincidence) also 
the 
characteristic operator function of $T^2$. This concludes the proof of the 
lemma.
\end{proof}

Note that the preceding result also shows that $T^2$ is a $C_0(2)$ operator.

\section{}

\indent

Our next step is to characterize in terms of $\Theta(\lambda)$ the 
reducibility 
of $T^2$.

\begin{lem}\label{lem2}
The operator $T^2$ is reducible if and only if there exist $Q_i = Q^*_i = 
Q^2_i$, $Q_i\in {\cl L}({\bb C}^2)$ $(i=1,2)$ so that
\begin{equation}\label{eq2.1}
\Theta(\lambda) Q_2 = Q_1 \Theta (\lambda) \qquad (\lambda\in {\bb D})
\end{equation}
and $0 \ne Q_i \ne I_{{\bb C}^2}$ $(i=1,2)$.
\end{lem}

\begin{proof}
If $Q_1,Q_2$ as above exist, then (since rank $Q_1 = 1 =$ rank $Q_2$) there 
exist unitary operators in ${\cl L}({\bb C}^2)$ so that
\begin{equation}\label{eq2.2}
W_1\Theta(\lambda)W_2 = \begin{bmatrix} \theta_1(\lambda)&0\\ 
0&\theta_2(\lambda)\end{bmatrix}\qquad (\lambda\in {\bb D}) \text{ for 
functions } \theta_1(\cdot), \theta_2(\cdot).
\end{equation}
Indeed, if $W_1$ and $W_2$ are unitary operators in ${\cl L}({\bb C}^2)$ such 
that
\[
Q_1{\bb C}^2 = W^*_1({\bb C}\oplus \{0\}), \quad Q_2{\bb C}^2 = W_2({\bb C} 
\oplus \{0\}),
\]
then
\begin{align*}
&W_1 \Theta(\lambda)W_2 \begin{bmatrix} 1&0\\ 0&0\end{bmatrix} - 
\begin{bmatrix} 
1&0\\ 0&0\end{bmatrix} W_1\Theta(\lambda) W_2=\\
=~ &W_1\left\{\Theta(\lambda)W_2 \begin{bmatrix} 1&0\\ 0&0\end{bmatrix} W^*_2 
- 
W^*_1 \begin{bmatrix} 1&0\\ 0&0\end{bmatrix} W_1 \Theta(\lambda)\right\} W_2 
=\\
=~ &W_1 (\Theta(\lambda) Q_2-Q_1\Theta(\lambda)) = 0.
\end{align*}
Thus ${\bb C}\oplus \{0\}$ (and hence also $\{0\}\oplus {\bb C}$) reduces 
$W_1\Theta(\lambda)W_2$ and consequently this operator has the form 
\eqref{eq2.2}.

Clearly the $\theta_1,\theta_2$ in \eqref{eq2.2} are inner (and non-constant). 
Let
\begin{equation}\label{eq2.3}
T_i = P_{H_i} S|_{\ds H_i}, \text{ where } H_i = H^2\ominus \theta_iH^2\qquad 
(i=1,2).
\end{equation}
Then the characteristic operator function of $T_1\oplus T_2$ is the right hand 
side of 
\eqref{eq2.2} which coincides with $\Theta(\lambda)$. Thus $T^2$ and 
$T_1\oplus 
T_2$ are unitarily equivalent.

Conversely, if $T^2$ is reducible then $T^2$ is unitarily equivalent to the 
direct sum $T'_1 \oplus T'_2$, where $T'_i = T^2|_{\textstyle H_i}$ $(i=1,2)$,  
$H_1,H_2$ 
are 
reducing subspaces for $T^2$, and $H=H_1\oplus H_2$. Clearly each $T'_i\in 
C_{00}$ and since 
the  defect indices of the $T'_i$s sum up to 2, it follows that each $T'_i\in 
C_0(1)$. Thus the characteristic operator function of $T'_1\oplus T'_2$ 
coincides with
\begin{equation}\label{eq2.4}
\begin{bmatrix} \theta_1(\lambda)&0\\ 0&\theta_2(\lambda)\end{bmatrix},
\end{equation}
where $\theta_i$ is the characteristic function of $T'_i$ $(i=1,2)$. Again 
$\Theta(\lambda)$ is connected to \eqref{eq2.4} by a relation of the form 
\eqref{eq2.2}, that is,
\[
\Theta(\lambda) \equiv W^*_1\begin{bmatrix} \theta_1(\lambda)&0\\ 
0&\theta_2(\lambda)\end{bmatrix} W^*_2,
\]
where $W_1,W_2$ are again unitary. Then
\[
Q_1 = W^*_1 \begin{bmatrix} 1&0\\ 0&0\end{bmatrix}W_1,\quad Q_2 = W_2 
\begin{bmatrix} 1&0\\ 0&0\end{bmatrix} W^*_2
\]
satisfy \eqref{eq2.1}
\end{proof}

\begin{rem}
Note that in \eqref{eq2.1}, the orthogonal projections $Q_1,Q_2$ are of rank 
one. Such a projection $Q$ is of the form
\begin{equation}\label{eq2.5}
Q = f\otimes f = \begin{bmatrix} |f_1|^2&f_1\bar f_2\\ f_2\bar f_1&|f_2|^2 
\end{bmatrix},
\end{equation}
where
\[
f = f_1\oplus f_2\in {\bb C}^2,\quad \|f\|=1.
\]
Thus
\begin{equation}\label{eq2.6}
Q = \begin{bmatrix} q&r\bar\theta\\ r\theta&1-q\end{bmatrix}, \text{ where } 
0\le q \le 1, |\theta| = 1, r=(q(1-q))^{1/2}.
\end{equation}
\end{rem}

\section{}

\indent

In this paragraph we study the relation \eqref{eq2.1} using the representation 
\eqref{eq2.6} for $Q=Q_i$ $(i=1,2)$ and the form \eqref{eq1.2} of 
$\Theta(\lambda)$. Thus we have
\begin{equation}\label{eq3.1}
\begin{bmatrix} b(\lambda)&\lambda d(\lambda)\\ d(\lambda)&b(\lambda) 
\end{bmatrix} \begin{bmatrix} q_2&r_2\bar\theta_2\\ 
r_2\theta_2&1-q_2\end{bmatrix} = \begin{bmatrix} q_1&r_1\bar\theta_1\\ 
r_1\theta_1&1-q_1\end{bmatrix} \begin{bmatrix} b(\lambda)&\lambda d(\lambda)\\ 
d(\lambda)&b(\lambda)\end{bmatrix},
\end{equation}
where
\begin{equation}\label{eq3.2}
0\le q_1,q_2 \le 1, |\theta_1| = |\theta_2| = 1, r_i = (q_i(1-q_i))^{1/2}\quad 
(i=1,2).
\end{equation}

We begin by noting that
\begin{equation}\label{eq3.3}
|b(\lambda)|^2 + |d(\lambda)|^2 \not\equiv 0\qquad (\lambda\in {\bb D}),
\end{equation}
since otherwise we would have $m(\lambda) \equiv 0$. In discussing 
\eqref{eq3.1} 
we will consider several cases:\medskip

\n {\bf Case I.} If $b(\lambda)\equiv 0$ $(\lambda\in {\bb D})$, then 
\eqref{eq3.1} becomes:
\[
\begin{bmatrix} \lambda d(\lambda)r_2\theta_2&\lambda d(\lambda) (1-q_2)\\
d(\lambda)q_2&d(\lambda)r_2\bar\theta_2\end{bmatrix} = \begin{bmatrix} 
r_1\bar\theta_1d(\lambda)&\lambda q_1d(\lambda)\\ 
(1-q_1)d(\lambda)&r_1\theta_1 
\lambda d(\lambda)\end{bmatrix}
\]
which is possible if and only if $r_1=0=r_2$ and $q_1= 1-q_2$. In this case 
$T^2$ is 
reducible.\medskip

\n {\bf Case II.} If $d(\lambda)\equiv 0$ $(\lambda\in {\bb D})$, then
\[
Q_2=Q_1 = \text{ any } Q = Q^* = Q^2 \text{ with rank } Q=1
\]
and again $T^2$ is reducible.\medskip

\n {\bf Case III.} If $b(\lambda) \not\equiv 0, d(\lambda)\not\equiv 0$ 
$(\lambda\in {\bb D})$, then \eqref{eq3.1} is equivalent to the equations
\begin{alignat*}{2}
b(q_2-q_1) &= d(r_1\bar\theta_1-\lambda r_2\theta_2),&\quad 
b(r_2\bar\theta_2-r_1\bar\theta_1) &= \lambda d(q_1+q_2-1)\\
d(q_1+q_2-1) &= b(r_1\theta_1-r_2\theta_2), &\quad b(q_2-q_1) &= 
d(r_2\bar\theta_2 - \lambda r_1\theta_1),
\end{alignat*}
which in turn are equivalent to
\begin{equation}\label{eq3.4}
\left\{\begin{array}{l}
r_1\theta_1 = r_2\theta_2,\quad  q_2+q_1=1\\
b(\lambda)(1-2q_1) \equiv d(\lambda) (\bar\theta_1-\lambda\theta_1) r_1 \quad 
(\lambda\in  {\bb D}).\end{array}\right.
\end{equation}
In \eqref{eq3.4}, $q_1=1/2$, if and only if $r_1=0$, i.e.\ $q_1=0$ or 1, a 
contradiction. Thus we can divide by $1-2q_1$ and \eqref{eq3.4} implies (with 
$\theta=\theta_1$)
\begin{equation}\label{eq3.5}
\left\{\begin{array}{l}
b(\lambda)\equiv d(\lambda) (\bar\theta-\lambda\theta)\rho\qquad (\lambda\in 
{\bb D})\\
\text{for some}\quad \rho\in {\bb R},\rho\ne 0.\end{array}\right.
\end{equation}

Conversely, if \eqref{eq3.5} holds, then setting
\[
q_1 = \frac12 \pm \frac12 \frac1{(4\rho^2+1)^{1/2}} \text{ (according to 
whether 
} \rho \lessgtr 0),
\]
and $q_2=1-q_1$, $\theta_2=\theta_1 =\theta$, we obtain \eqref{eq3.4}.

We now summarize our discussion in terms of $m(\cdot)$ (see \eqref{eq1.3a}, 
\eqref{eq1.3b}), instead of $b(\cdot)$ and $d(\cdot)$, obtaining the 
following:

\begin{lem}\label{lem3}
The operator $T^2$ is reducible if and only if one of the following conditions 
holds:
\begin{align}
\label{eq3.6}
m(-\lambda) &\equiv -m(\lambda)\quad (\forall\lambda\in {\bb D})\qquad 
\text{(Case I above);}\\
\label{eq3.7}
m(-\lambda) &\equiv m(\lambda)\quad (\forall \lambda\in {\bb D})\qquad 
\text{(Case II above);}
\end{align}
or there exist $\rho\in {\bb R}$, $\rho\ne 0$ and $\theta\in {\bb C}$, 
$|\theta|=1$, such that the function
\begin{subequations}
\begin{align}
\label{eq3.8a}
n(\lambda) &\equiv m(\lambda) (\rho\theta \lambda^2 + \lambda-\rho \bar\theta) 
\qquad (\lambda\in{\bb D})\\
\intertext{satisfies}
\label{eq3.8b}
n(\lambda) &\equiv n(-\lambda)\qquad (\lambda\in {\bb D})\qquad \text{(Case 
III 
above).}
\end{align}
\end{subequations}
\end{lem}

\section{}

\indent

We shall now give a more transparent form to  conditions \eqref{eq3.8a}, 
\eqref{eq3.8b} above. To this end note that
\[
\rho\theta\lambda^2 + \lambda-\rho\bar\theta \equiv \rho\theta(\lambda - 
\delta_+ \bar\theta)(\lambda-\delta_-\bar\theta),
\]
where
\begin{equation}\label{eq4.1}
\delta_\pm = \frac{-1\pm\sqrt{4\rho^2+1}}{2\rho}.
\end{equation}
Thus (with $\mu = \bar\theta\delta_+$), we have
\begin{equation}\label{eq4.2}
\rho\theta\lambda^2 + \lambda-\rho\bar\theta = 
-\rho\delta_-(\lambda-\mu)(1+\bar\mu\lambda).
\end{equation}

Using this representation in\eqref{eq3.8a},  condition \eqref{eq3.8b} becomes
\[
m(\lambda)(\lambda-\mu)(1+\bar\mu\lambda) \equiv m(-\lambda)(-\lambda-\mu) 
(1-\lambda\bar\mu)\qquad (\lambda\in {\bb D}),
\]
which can be written  (since $0<|\mu|<1$) as
\begin{equation}\label{eq4.3}
m(\lambda) \frac{\lambda-\mu}{1-\bar\mu\lambda} \equiv m(-\lambda) 
\frac{(-\lambda)-\mu}{1-\bar\mu(-\lambda)}\qquad (\lambda\in {\bb D}).
\end{equation}
Thus $m(-\mu) = 0$ and therefore
\begin{equation}\label{eq4.4}
m(\lambda) = p(\lambda) \frac{\lambda+\mu}{1+\bar\mu\lambda}\qquad (\lambda\in 
{\bb D}),
\end{equation}
where $p(\cdot)\in H^\infty$ is an (other) inner function. Obviously 
\eqref{eq4.3} is equivalent to
\begin{equation}\label{eq4.5}
p(\lambda)\equiv p(-\lambda)\qquad (\lambda\in {\bb D}).
\end{equation}
This discussion together with Lemma \ref{lem3}, readily yields the following

\begin{thm}\label{thm1}
The operator $T^2$ is reducible iff either
\begin{equation}\label{eq4.6}
m(\lambda) = m(-\lambda)\qquad (\lambda\in {\bb D})
\end{equation}
or there exists a $\mu\in {\bb D}$ such that
\begin{equation}\label{eq4.7}
m(\lambda)\equiv p(\lambda) \frac{\lambda+\mu}{1+\bar\mu\lambda} \qquad 
(\lambda\in {\bb D}),
\end{equation}
where $p(\cdot)\in H^\infty$ satisfies
\begin{equation}\label{eq4.8}
p(\lambda)\equiv p(-\lambda)\qquad (\lambda\in {\bb D}).
\end{equation}
\end{thm}

\begin{rem}
Case \eqref{eq3.6} is contained in the second alternative above when 
$\mu=0$.
\end{rem}

\section{}

\indent

In order to study the lattices $\text{Lat}\{T^2\}$ and $\text{Lat}\{T^2\}'$ we 
first bring together the following characterization of the $C_0(N)$ operators 
that are multiplicity free.

\begin{pro}\label{pro1}
Let $\widetilde T$ be a $C_0(N)$ operator. Then the following statements are 
equivalent.
\begin{itemize}
\item[\rm (1)] $\widetilde T$ is multiplicity free (that is, $\widetilde T$ has 
a 
cyclic vector).
\item[\rm (2)] $\text{\rm Lat}\{\widetilde T\} = \text{\rm Lat}\{\widetilde 
T\}'$.
\item[\rm (3)] The minors of  the characteristic matrix function of order $N-1$  
have 
no common inner divisor.
\end{itemize}
\end{pro}

\begin{proof}
The equivalence of (1) and (3) is contained in the equivalence of (i) and (ii) 
in Theorem \ref{thm2} in \cite{S-F2}. The implication (1) implies (2) is an easy 
corollary of the implication (i) implies (vi) of the same theorem and is 
contained in Corollary 2.14 in Chapter 3 of \cite{B}. Finally, implication (3) 
implies (1) proceeds from the following lemma.
\end{proof}

\begin{lem}\label{lem4}
Let $T$ be an $C_0$ operator on the Hilbert space ${\cl H}$ and $f$ a maximal 
vector for $T$. Then $f$ is cyclic for $\{T\}'$.
\end{lem}

\begin{proof}
Let ${\cl M}$ be the cyclic subspace for $\{T\}'$ generated by $f$ and write 
$T\sim \left(\begin{smallmatrix} T'&X\\ 0&T''\end{smallmatrix}\right)$ for the 
decomposition ${\cl H} = {\cl M} \oplus {\cl M}^\bot$. Since ${\cl M}$ is 
hyperinvariant for $T$, it follows from Corollary 2.15 in Chapter 4 of \cite{B}, 
that the minimal functions satisfy $m_T = m_{T'} \cdot m_{T''}$. However, 
$f$ maximal for $T$ implies that $m_{T'} = m_T$ and hence $m_{T''} =1$. 
Therefore, ${\cl M}^\bot = (0)$ or ${\cl M} = {\cl H}$ which completes the 
proof.
\end{proof}

\section{}

\indent

Our next aim is to characterize the case when the operator $T^2$ is 
multiplicity 
free. According to Proposition \ref{pro1} that happens if and only if
\[
b(\lambda), d(\lambda) \text{ and } \lambda d(\lambda)
\]
have no common nontrivial inner divisor. Let $q(\lambda)$ be an inner divisor 
of 
$b(\lambda)$ and $d(\lambda)$, that is,
\begin{subequations}
\begin{align}
\label{eq6.1a}
m(\sqrt\lambda) + m(-\sqrt\lambda) &\equiv q(\lambda) r(\lambda)\\
\label{eq6.1b}
m(\sqrt\lambda) - m(-\sqrt\lambda) &\equiv q(\lambda)\lambda s(\lambda)
\qquad \raisebox{3ex}{$(\lambda\in{\bb D})$}
\end{align}
\end{subequations}
for some $r,s\in H^\infty$. It follows that
\begin{equation}\label{eq6.2}
m(\lambda)\equiv q(\lambda^2) (r(\lambda^2) -\lambda s(\lambda^2)),
\end{equation}
that is,  $m(\lambda)$ has an even inner divisor.

Conversely, if $m(\cdot)$ has an inner divisor (in $H^\infty$) $p(\cdot)$ 
satisfying
\begin{equation}\label{eq6.3}
p(\lambda)\equiv p(-\lambda),
\end{equation}
then $q(\lambda) = p(\sqrt\lambda) = p(-\sqrt\lambda)$ is in $H^\infty$ and 
inner. Thus $m(\lambda)$ can be represented as in \eqref{eq6.2} and clearly 
\eqref{eq6.2} implies \eqref{eq6.1a}, \eqref{eq6.1b}. Thus we obtained the 
following:

\begin{thm}\label{thm2}
The operator $T^2$ is  multiplicity free iff  the characteristic function 
$m(\lambda)$ for $T$ has no nontrivial inner divisor $p(\lambda)$ in $H^\infty$ 
such that (see \eqref{eq6.3})
\[
p(\lambda)\equiv p(-\lambda)\qquad (\forall\lambda\in {\bb D}).
\]
\end{thm}

\section{}

\indent

Our main result is now a direct consequence of Theorems \ref{thm1} and 
\ref{thm2} 
and Proposition \ref{pro1}, namely

\begin{thm}\label{thm3}
Let $T\in C_0(1)$ satisfy:
\begin{equation}\label{eqA}
m_T(\lambda)\not\equiv m_T(-\lambda)\tag{\rm A}
\end{equation}
{\rm (B)}~~For $m_T(\lambda_0)=0$, $\lambda_0\in {\bb D}$, the function
\[
m_{T,\lambda_0}(\lambda)  = m_T(\lambda)\Big/ \frac{\lambda-\lambda_0}{1- 
\bar\lambda_0\lambda}\qquad (\lambda\in {\bb D})
\]
is not even, that is,
\[
m_{T,\lambda_0}(\lambda)\not\equiv m_{T,\lambda_0}(-\lambda).
\]
{\rm (C)}~~There exists a nontrivial inner divisor $p(\lambda)$ (in 
$H^\infty$) 
of 
$m_T(\lambda)$ such that
\[
p(\lambda) \equiv p(-\lambda).
\]
Then
\begin{align}
&T^2 \text{ is irreducible, and}\tag{\rm D}\\
&\text{\rm Lat } T^2\ne \text{\rm Lat}\{T^2\}'.\tag{\rm E}
\end{align}
\end{thm}

\section{}

\n {\bf Remarks}

a)~~Let
\begin{equation}\label{eq8.1}
m_T(\lambda) = \frac{\lambda^2-\lambda_1}{1-\bar\lambda_1\lambda^2} 
\left(\frac{\lambda-\lambda_2}{1-\bar\lambda_2\lambda}\right)^2\qquad 
(\lambda\in {\bb D}),
\end{equation}
where $\lambda_1,\lambda_2\in {\bb D}$, $\lambda^2_2\ne \lambda_1$. Then $m$ 
fulfills the solutions (A), (B), (C) in Theorem \ref{thm3},
  $T^2$ satisfies (D) and (E)  above and hence $\dim H = 4$.

b)~~If $\dim H=3$ then
\[
m_T(\lambda) = \frac{\lambda-\lambda_1}{1-\bar\lambda_1\lambda} 
\frac{\lambda-\lambda_2}{1-\bar\lambda_0\lambda} 
\frac{\lambda-\lambda_3}{1-\bar 
\lambda_3\lambda}\qquad (\lambda\in {\bb D})
\]
with some $\lambda_1,\lambda_2,\lambda_3\in {\bb D}$. If $m_T$ satisfies (C) 
then $(\lambda_1+\lambda_2)(\lambda_2+\lambda_3)(\lambda_3+\lambda_1)=0$ and 
$m_T(\lambda)$ has the form (upon relabelling the $\lambda_i$'s)
\begin{equation}\label{eq8.2}
m_T(\lambda) =\frac{\lambda^2-\lambda^2_1}{1-\bar\lambda^2_1\lambda^2} 
\frac{\lambda-\lambda_2}{1-\bar\lambda_2\lambda}\qquad (\lambda\in {\bb D}).
\end{equation}
Consequently $m_T$ does not satisfy (B).  Thus for Theorem \ref{thm3} to 
hold 
it is necessary that $\dim H\ge 4$.

3)~~Let $m_T$ be singular, that is,
\[
m_T(\lambda) = \exp\left[-\frac1{2\pi} \int\limits^\pi_0 
\frac{e^{it}+\lambda}{e^{it}-\lambda} d\mu(e^{it})\right]
\]
with $\mu$ a singular measure on $\partial{\bb D} = \{e^{it}\colon\ 0\le t < 
2\pi\}$. Assume that there exists a Borel set $\Omega\subset \partial{\bb D}$ 
so 
that
\[
\mu(\Omega) = \mu(\partial{\bb D}),\quad \mu(\{\bar\lambda\colon \ \lambda\in 
\Omega\}) = 0.
\]
(e.g.\ $\mu = \delta_1$, the point mass at 1). Then
\begin{equation}\label{eq8.3}
\text{Lat}\{T^2\} = \text{Lat}\{T^2\}' = \text{Lat}\{T\}.
\end{equation}
Indeed, in this case (C) above does not hold.

\vspace{1in}

\n Department of Mathematics\\
Texas A\&M University\\
rdouglas@math.tamu.edu
\end{document}